\documentclass{ifacconf}

\usepackage{graphicx}      
\usepackage{natbib}        
\usepackage{amsmath} 
\usepackage{amsfonts}
\usepackage{amssymb}
\usepackage{enumerate}
\usepackage{graphicx}
\usepackage{epstopdf}
\usepackage{subfigure}
\usepackage{url} 
\usepackage{tikz}
\usepackage{color}
\usepackage{float}

\newtheorem{theorem}{Theorem}

\newtheorem{definition}{Definition}
\begin{document}
\begin{frontmatter}

\title{Event-Triggered Newton-Based Extremum Seeking Control} 


\author[First]{Victor Hugo Pereira Rodrigues}
\author[First]{Tiago Roux Oliveira} 
\author[Second]{Miroslav Krsti{\' c}} 
\author[Third]{Paulo Tabuada} 

\address[First]{Department of Electronics and Telecommunication Engineering,\\ State University of Rio de Janeiro (UERJ), Rio de Janeiro--RJ, Brazil. \\ \mbox{(e-mail: victor.rodrigues@uerj.br, tiagoroux@uerj.br)}}

\address[Second]{Department of Mechanical and Aerospace  Engineering,\\ University of California at San Diego, La Jolla--CA, USA. \\ \mbox{(e-mail:  mkrstic@ucsd.edu)}}

\address[Third]{Department of Electrical and Computer Engineering,\\ University of California at Los Angeles, Los Angeles--CA, USA. \\ \mbox{(e-mail:  tabuada@ee.ucla.edu)}}

\begin{abstract}                
This paper proposes the incorporation of static event-triggered control in the actuation path of Newton-based extremum seeking and its comparison with the earlier gradient version. As in the continuous methods, the convergence rate of the gradient approach depends on the unknown Hessian of the nonlinear map to be optimized, whereas the proposed event-triggered Newton-based extremum seeking eliminates this dependence, becoming user-assignable. This is achieved by means of a dynamic estimator for the Hessian's inverse, implemented as a Riccati equation filter. Lyapunov stability and averaging theory for discontinuous systems are applied to analyze the closed-loop system. Local exponential practical stability is guaranteed to a small neighborhood of the extremum point of scalar and static maps.  Numerical simulations illustrate the advantages of the proposed approach over the previous gradient method, including improved convergence speed, followed by a reduction in the amplitude and updating frequency of the control signals.
\end{abstract}

\begin{keyword}
Extremum Seeking, Event-Triggered Control, Averaging Theory, Scalar Maps, Exponential Stability.
\end{keyword}

\end{frontmatter}

\section{Introduction}
In 1922, the concept of extremum seeking was introduced by the French engineer Maurice Leblanc to maintain the maximum efficiency of power transfer from an electrical transmission line to a tram car \citep{L:1922}. Almost eighty years later, a rigorous  stability analysis of extremum seeking feedback was proposed in \citep{KW:2000}, using 
averaging and singular perturbation methods. This method optimizes the performance of an unknown dynamic system by employing probing and demodulation periodic signals in order to estimate the gradient and search for the extremum point of the corresponding output cost function \citep{IJACSP2016}.

After that, numerous studies on the topic have been published over the time, including theoretical advances and engineering practice
\citep{AK:2003,zhang2012extremum,Aminde2013,c23,AJC2014}.  In the field of applications, extremum seeking has been used in autonomous vehicles and mobile sensors \citep{SD:2009} and optimization of PID controllers applied to Neuromuscular Electrical Stimulation (NMES) \citep{TRO:2019,TRO:2020}. Recent works have also employed extremum seeking to derive policies for Nash equilibrium seeking in non-cooperative games \citep{FKB:2012}, including cases where the players actions are subject to distinct dynamics governed by diffusive Partial Differential Equations (PDEs) \citep{TRoux:2021may} and with delays \citep{TRoux:2021dec}. As illustrated in the book \citep{TRoux:2022}, this approach could be extended to other classes of PDEs. 

On the other hand, due to the rapid development of network technologies, communication networks between devices and plants have become fast and reliable, reducing wiring costs and simplifying maintenance \citep{ZHGDDYP:2020}. However, one significant disadvantage remains: network bandwidth is limited, and network traffic congestion is unavoidable, leading to delays and packet loss \citep{HNX:2007}. To mitigate this problem, event-triggered Control (ETC) can be employed \citep{HJT:2012}. In contrast to periodic sample-date controllers, ETC execute control tasks only when a predefined condition based on the plant’s state is triggered, thereby reducing the consumption of computational resources with guaranteed asymptotic stability properties \citep{T:2007}.

In this context, combining extremum seeking with event-triggered control maybe interesting in order to inherit the advantages of both methods. Oriented for static maps, an event-triggered control scheme for scalar extremum seeking was proposed in \citep{VHPR:2022} and \citep{VHPR:2023b}, where the closed-loop stability of the closed-loop system was guaranteed without Zeno behavior. Later on, both static and dynamic triggering scenarios were developed \citep{VHPR:2023a} showing the benefits from this combination for the multivariable case.

The references cited above \citep{VHPR:2022,VHPR:2023b,VHPR:2023a} employed gradient-based extremum seeking to maximize or minimize the objective function, where the algorithm's convergence rate speed is proportional to the Hessian (second derivative) of the nonlinear map to be optimized. In 2010, for a single-input case, an estimate of the second derivative was employed in a Newton-like method \citep{MMB:2010} to reduce sensitivity to the curvature of the plant map. The Newton method requires the knowledge of the inverse of the Hessian, which is not a trivial task in model-free optimization. To mitigate this, \citep{GKN:2012} proposed a dynamic estimator in the form of a Riccati differential equation for the Hessian's inverse matrix, eliminating the dependence on the unknown second derivative in the the convergence rate, and making it user-assignable. The key difference between the gradient-based and Newton-based methods is that the former depends on the Hessian, while the latter is independent of it. This offers a significant advantage to the extremum seeking method, where the Hessian is unknown.

In this paper, we focus on the design and analysis of a scalar Newton-based extremum seeking method for static maps within an event-triggered framework, aiming to leverage the benefits of improved convergence rates independent of the Hessian, by employing a Riccati equation filter to estimate its inverse. Additionally, we examine the impact of control signal updates and the resulting input-output responses. To address this, a Lyapunov-based  criterion and averaging theory for discontinuous systems are utilized to characterize the stability of the closed-loop system. A numerical example illustrates the theoretical results.

\emph{Notation:} 
Throughout the paper, the absolute value of scalar variables \textcolor{black}{is} denoted by single bars $|\cdot|$. Consider $\varepsilon \in \lbrack -\varepsilon_{0}\,, \varepsilon_{0} \rbrack \subset \mathbb{R}$ and the mappings $\delta_{1}(\varepsilon)$ and $\delta_{2}(\varepsilon)$, where $\delta_{1}: \lbrack -\varepsilon_{0}\,, \varepsilon_{0} \rbrack \to \mathbb{R}$ and $\delta_{2}: \lbrack -\varepsilon_{0}\,, \varepsilon_{0} \rbrack \to \mathbb{R}$. \textcolor{black}{The function $\delta_1(\varepsilon)$ has magnitude of order $\delta_2(\varepsilon)$, denoted by} $\delta_{1}(\varepsilon) = \mathcal{O}(\delta_{2}(\varepsilon))$, if there exist positive constants $k$ and $c$ such that $|\delta_{1}(\varepsilon)| \leq k |\delta_{2}(\varepsilon)|$, for all $|\varepsilon|<c$ \citep{K:2002,IJRNC2011}.

\section{Problem Formulation} \label{sec:prblFrm_siso}

\begin{figure}[!ht]
	\centering
	\includegraphics[width=8.5cm]{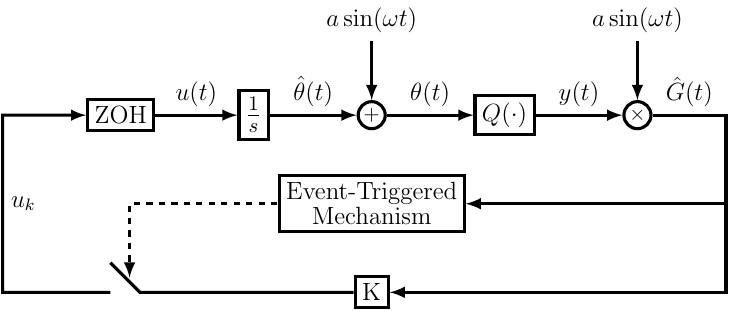}
	\caption{Static Event-Triggered Gradient-based Extremum Seeking \citep{VHPR:2022}.}
\label{fig:BD_SET_GradientES_SISO}
\end{figure}

    Fig.~\ref{fig:BD_SET_GradientES_SISO} shows the complete block diagram of the Static Event-Triggered Gradient-based Extremum Seeking (SET-GradientES) proposed in \citep{VHPR:2022}, where the algorithm's convergence rate can be highly sensitive to the curvature \textcolor{black}{of the map $Q(\cdot)$ to be optimized.} 
    The SET-GradientES may require conservative tuning to ensure stability conditions within a given operating domain, resulting in slow optimization responses. The SET-GradientES algorithm is locally convergent with  convergence rate determined by the unknown Hessian matrix of the static map.

\begin{figure}[!ht]
	\centering
	\includegraphics[width=8.5cm]{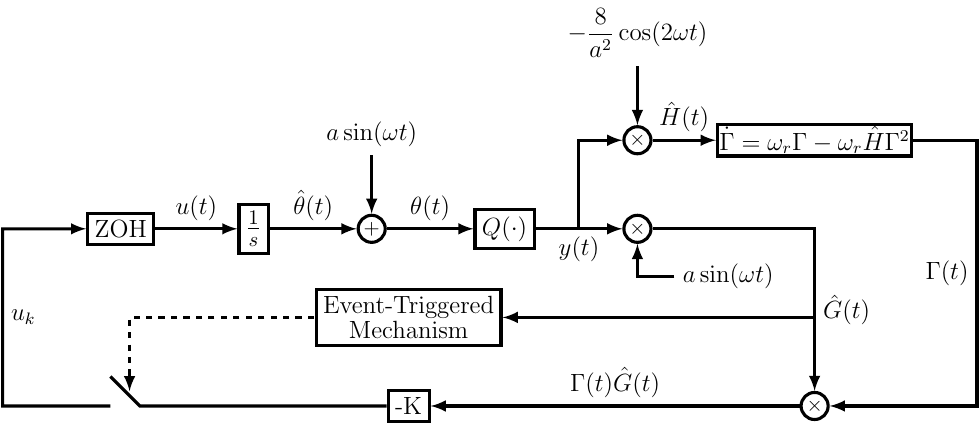}
\caption{Static Event-Triggered Newton-based Extremum Seeking.}
\label{fig:BD_SET_NewtonES_SISO}
\end{figure}

To overcome this limitation in SET-GradientES, here we propose the Static Event-Triggered Newton-based\linebreak Extremum Seeking (SET-NewtonES) depicted in Fig.~\ref{fig:BD_SET_NewtonES_SISO}. The SET-NewtonES executes the control task aperiodically in response to a triggered condition designed as a function of the gradient estimate while ensures that the convergence rate is arbitrarily specified by the user and remains unaffected by the Hessian $H^{\ast}$ of the map, which is unknown. 

We define \textcolor{black}{the output of the nonlinear static map $Q(\theta(t))$ as
\begin{align}
y(t)=Q(\theta(t)) = Q^{\ast}+\frac{H^{\ast}}{2}(\theta(t)-\theta^{\ast})^{2}\,, \label{eq:y_event_siso} 
\end{align}
}
$\!\!\!$where \textcolor{black}{$Q^{\ast} \in \mathbb{R}$ is a constant}, $H^{\ast} \in \mathbb{R}-\{0\}$ is the Hessian, $\theta^{\ast} \in \mathbb{R}$ is the unknown optimizer, and $\theta(t)\in \mathbb{R}$ is the input of the map. 

\textcolor{black}{ Although (\ref{eq:y_event_siso}) is a polynomial in $\theta$, it cannot be identified from a finite number of input/output pairs $(\theta, y) $ due to the real-time nature of extremum seeking optimization. Furthermore, extremum seeking is a versatile approach that can handle any unknown analytic function $Q(\theta)$, which can be expressed as a Taylor series expansion around a point $\theta^{\ast}$, where the function attains a minimum or maximum. This assumption enables a local quadratic approximation of the nonlinear map $ Q(\theta)$, serving as a foundational principle and justifying our focus on quadratic maps. Therefore, while the analysis was conducted for a quadratic map, our approach is not limited to quadratic functions $Q(\theta)$.}

The parameters of the nonlinear static map (\ref{eq:y_event_siso}) are not explicitly known; however, we can measure the output $y(t)$ and design the input signal $\theta(t)$. The input $\theta(t)$ is obtained from the real-time estimate $\hat{\theta}(t) \in \mathbb{R}$ of $\theta^{\ast}$ additively perturbed by the periodic signal $a \sin(\omega t)$, {\it i.e.},
\begin{align}
\theta(t)=\hat{\theta}(t)+a \sin(\omega t)\,. \label{eq:theta_event_siso}
\end{align}


\subsection{Continuous-Time Extremum Seeking}

Let us define the estimation error 
\begin{align}
\tilde{\theta}(t)=\hat{\theta}(t)-\theta^{\ast}\,, \label{eq:thetaTilde_event_siso}
\end{align}
and the gradient estimate
\begin{align}
\hat{G}(t)=a \sin(\omega t) ~y(t)\,, \label{eq:hatG_event_siso}
\end{align}
\textcolor{black}{by the demodulation signal}, $a \sin(\omega t)$, which has nonzero amplitude $a$ and frequency $\omega$ \citep{GKN:2012,K:2014}.

From (\ref{eq:theta_event_siso}) and (\ref{eq:thetaTilde_event_siso}), we can write
\begin{align}
\theta(t)=\tilde{\theta}(t)+a \sin(\omega t)+\theta^{\ast}\,, \label{eq:theta_2_event_siso}
\end{align}
and, therefore, by plugging (\ref{eq:theta_2_event_siso}) into (\ref{eq:y_event_siso}), $y(t)$ can also be written as
\begin{align}
y(t)
&=Q^{\ast}+\frac{H^{\ast}a^{2}}{4}+\frac{H^{\ast}}{2}\tilde{\theta}^{2}(t)+a \sin(\omega t)H^{\ast}\tilde{\theta}(t)\nonumber\\
&\quad-\frac{H^{\ast}a^{2}}{4} \cos(2\omega t)\,.
 \label{eq:y_1_event_siso}
\end{align}
Thus, \textcolor{black}{from (\ref{eq:hatG_event_siso}) and (\ref{eq:y_1_event_siso}),} the gradient estimate is \textcolor{black}{given by}
\begin{align}
\hat{G}(t)&= \frac{a^2 H^{\ast}}{2}\left(1-\cos\left(2\omega t\right)\right)\tilde{\theta}(t)+\frac{a H^{\ast}}{2}\sin\left(\omega t\right)\tilde{\theta}^{2}(t)\nonumber \\
&\quad+\left(a Q^{\ast}+\frac{3 a^{3}H^{\ast}}{8}\right)\sin\left(\omega t\right)-\frac{a^{3}H^{\ast}}{8}\sin\left(3\omega t\right)\,. \label{eq:hatG_2_event_siso}
\end{align}
Notice the quadratic term in $\tilde{\theta}(t)$ of (\ref{eq:hatG_2_event_siso}) \textcolor{black}{may be neglected} in a local analysis \citep{AK:2003}. Thus, hereafter the gradient estimate \textcolor{black}{is given by}
\begin{align}
\hat{G}(t)&=\frac{a^2 H^{\ast}}{2}\left(1-\cos\left(2\omega t\right)\right)\tilde{\theta}(t)\nonumber \\
&\quad+\left(a Q^{\ast}+\frac{3 a^{3}H^{\ast}}{8}\right)\sin\left(\omega t\right)-\frac{a^{3}H^{\ast}}{8}\sin\left(3\omega t\right) \,. \label{eq:hatG_4_event_siso}
\end{align}

On the other hand, from the time-derivative of (\ref{eq:thetaTilde_event_siso}) and the SET-NewtonES scheme of Fig.~\ref{fig:BD_SET_NewtonES_SISO}, the dynamics governing  $\hat{\theta}(t)$, as well as $\tilde{\theta}(t)$, is given by
\begin{align}
\dot{\tilde{\theta}}(t)&=\dot{\hat{\theta}}(t)=u(t) \label{eq:dotTildeTheta_1_event_siso}\,, 
\end{align}
where $u$ is the SET-NewtonES law to be designed as
\begin{align}
u(t)=-K\Gamma(t)\hat{G}(t)\,, \quad K>0\,, \quad \forall t\geq 0\,.
\end{align}

By computing the time-derivative of (\ref{eq:hatG_4_event_siso}) and considering (\ref{eq:dotTildeTheta_1_event_siso}), one gets
\begin{align}
\dot{\hat{G}}(t)
&=\frac{a^{2}H^{\ast}}{2}\left(1-\cos\left(2\omega t\right)\right)u(t)+a^{2}\omega H^{\ast}\sin\left(2\omega t\right)\tilde{\theta}(t) \nonumber \\
&\quad+\left(\! a\omega Q^{\ast}\!+\!\frac{3a^{2}\omega H^{\ast}}{8}\!\right)\!\cos\left(\omega t\right)\!-\!\frac{3a^{2}\omega H^{\ast}}{8}\!\cos\left(3 \omega t\right) \,. \label{eq:dotHatG_event_siso}
\end{align}

\subsection{Hessian Estimate and its Inverse}

The SET-NewtonES strategy relies on two key components; the \textit{Hessian estimate}: 
\begin{align}
&\hat{H}(t)=-\frac{8}{a^2}\cos(2\omega t)y(t) \nonumber \\
&=H^{\ast}(1+\cos(4\omega t)-2\cos(2\omega t))-\frac{8Q^{\ast}}{a^2}\cos(2\omega t)\nonumber \\
&\quad-\frac{4H^{\ast}}{a}(\sin(3\omega t)+\sin(\omega t))\tilde{\theta}(t)-\frac{4H^{\ast}}{a}\cos(2\omega t)\tilde{\theta}^2(t) \label{eqhatH}\,,
\end{align}
and the dynamic \textit{Riccati equation}:
\begin{align} 
\dot{\Gamma}\left(t\right) = \omega_{r}\Gamma\left(t\right) - \omega_{r}\hat{H}\left(t\right)\Gamma^2\left(t\right) \,, \quad \omega_{r}>0\,, \label{eq:Gamma} 
\end{align}
whose solution gives an estimate of the Hessian's inverse $H^{\ast-1}$, ensuring robustness even when the Hessian estimate $\hat{H}(t)$ becomes singular (crossing by zero). 

\section{Event-Triggered Control Emulation of Newton-Based Extremum Seeking}


\textcolor{black}{Let $t_{k}$ denote \textcolor{black}{an} unbounded monotonically increasing sequence of time, {\it i.e.}, }
\begin{align}
0=t_{0}<t_{1}<\ldots<t_{k}<\ldots\,, \quad k \in \mathbb{N}\,, \lim_{k \to \infty} t_{k}=\infty \,, \label{eq:s_k_event_siso}
\end{align}
with the aperiodic sampling intervals $\tau_{k}=t_{k+1}-t_{k}>0$. 

We consider continuous measurement of the system output while actuating the system using an event-based approach. The actuator transforms the discrete-time control input $u_k=U(t_{k})$ to a piece-wise continuous control input $u(t)$ as in sampled data systems with a zero-order hold mechanism. By assuming that there is no delay in the Sensor-to-Controller and Controller-to-Actuator branches, we define the control input
\begin{align}
u(t)\!=\!-K\Gamma(t_{k})\hat{G}(t_{k})\,, \quad K\!>\!0\,, \quad \forall t \!\in\! \lbrack t_{k}\,, t_{k+1}\phantom{(}\!\!) \,, \quad k\!\in\! \mathbb{N}  \,, \label{eq:U_event_siso}
\end{align}
and, since the control updates are made at discrete times rather than continuously, we introduce the error 
\begin{align}
e(t):=\Gamma(t_{k})\hat{G}(t_{k})-\Gamma(t)\hat{G}(t) \,, \quad \forall t \in \lbrack t_{k}\,, t_{k+1}\phantom{(}\!\!) \,, \quad k\in \mathbb{N} \,. \label{eq:e_event_siso}
\end{align}
Therefore, by using (\ref{eq:e_event_siso}), the SET-NewtonES control law (\ref{eq:U_event_siso}) can be rewritten as
\begin{align}
u(t)=-K\Gamma(t)\hat{G}(t)-Ke(t)\,, \quad \forall t \in \lbrack t_{k}\,, t_{k+1}\phantom{(}\!\!) \,, \quad k\in \mathbb{N}  \,. \label{eq:U_event_v2_siso}
\end{align}

Now, plugging (\ref{eq:U_event_v2_siso}) into (\ref{eq:dotTildeTheta_1_event_siso}) and (\ref{eq:dotHatG_event_siso}), we arrive at the following Input-to-State Stable (ISS) \citep{K:2002} representations for the dynamics of $\hat{G}(t)$ and $\tilde{\theta}$ with respect to the error vector $e(t)$ in (\ref{eq:e_event_siso}): 
\begin{align}
&\dot{\hat{G}}(t)\mathbb{=} -\frac{a^{2}}{2}\left(1-\cos\left(2\omega t\right)\right)H^{\ast}K\Gamma(t)\hat{G}(t)\nonumber \\
&-\frac{a^{2}}{2}\left(1-\cos\left(2\omega t\right)\right)H^{\ast}Ke(t)-\frac{3a^{2}\omega H^{\ast}}{8}\cos\left(3 \omega t\right) \nonumber \\
&+a^{2}\omega H^{\ast}\sin\left(2\omega t\right)\tilde{\theta}(t)+\left(a\omega Q^{\ast}+\frac{3a^{2}\omega H^{\ast}}{8}\right)\cos\left(\omega t\right)\,, \label{eq:dotHatGav_event_3_siso} \\
&\dot{\tilde{\theta}}(t)\mathbb{=}-K\Gamma(t)\hat{G}(t)-Ke(t) \nonumber \\
&\mathbb{=} -\frac{a^2 }{2}\left(1-\cos\left(2\omega t\right)\right)K\Gamma(t)H^{\ast}\tilde{\theta}(t)-Ke(t)\nonumber \\
&+\left[-\left(a Q^{\ast}+\frac{3 a^{3}H^{\ast}}{8}\right)\sin\left(\omega t\right)+\frac{a^{3}H^{\ast}}{8}\sin\left(3\omega t\right)\right]K\Gamma(t)\,. \label{eq:dotTildeTheta_2_event_siso}
\end{align}

In a conventional sampled-data implementation, the transmission times are distributed equidistantly in time, meaning that $t_{k+1}=t_{k}+ h$, for all \textcolor{black}{$k$} , and some interval $h>0$. In event-triggered control, however, these transmission are orchestrated by a \textcolor{black}{monitoring} mechanism that invokes transmissions when the difference between the current value of the output and its previously transmitted value becomes too large in an appropriate sense \citep{HJT:2012}. 
In Definition~\ref{def:staticEvent_siso} below, our triggering strategy is presented.

\begin{definition} \label{def:staticEvent_siso}
Consider the nonlinear mapping $\Xi: \mathbb{R} \times \mathbb{R}\mapsto \mathbb{R}$ given by 
\begin{align}
\Xi(\hat{G},e) = \sigma|\hat{G}(t)|-\beta|e(t)|\,, \quad \sigma \in (0,1)\,, \beta>0\,,   
\end{align}
with $K>0$ being the control gain in (\ref{eq:U_event_siso}). The event-triggered controller with triggering condition consists of two components:
\begin{enumerate}
	\item \textcolor{black}{A sequence of increasing times} $I=\{t_{0}\,, t_{1}\,, t_{2}\,,\ldots\}$, with $t_{0}=0$, generated under the following rules 
		\begin{itemize}
			\item If $\left\{t \in\mathbb{R}^{+}: t>t_{k} ~ \wedge ~ \Xi (\hat{G},e) < 0 \right\} = \emptyset$, then the set of the times of the events is $I=\{t_{0}\,, t_{1}\,, \ldots, t_{k}\}$.
			\item If $\left\{t \in\mathbb{R}^{+}: t>t_{k} ~ \wedge ~ \Xi (\hat{G},e) < 0  \right\} \neq \emptyset$, the next event time is given by
				\begin{align}
					t_{k+1}&=\textcolor{black}{\inf}\left\{t \in\mathbb{R}^{+}: t>t_{k} ~ \wedge ~ \Xi (\hat{G},e) <0 \right\}\,, \label{eq:tk+1_event} 
				\end{align}
				which is the event-trigger mechanism.
		\end{itemize}
	\item A feedback control action updated at the generated triggering instants given by
		\begin{align}
			u_{k}=-K\Gamma(t_{k})\hat{G}(t_{k}) \,,  \quad \forall t \in \lbrack t_{k}\,, t_{k+1}\phantom{(}\!\!)\,, \quad k\in \mathbb{N}\,. \label{eq:U_MD2}
		\end{align}
\end{enumerate}  
\end{definition}

\section{Time-Scaling and Closed-Loop\\ Average System}

By using the transformation $\bar{t}=\omega t$, where 
\begin{align}
\textcolor{black}{T}&\textcolor{black}{:=\frac{2\pi}{\omega}}\,, \label{eq:omega_event_1_siso}
\end{align} 
it is possible to rewrite the dynamics (\ref{eq:dotHatGav_event_3_siso})--(\ref{eq:dotTildeTheta_2_event_siso}) in a different time-scale such that
\textcolor{black}{ 
\begin{align}
\frac{d\hat{G}(\bar{t})}{d\bar{t}}&=\frac{1}{\omega}\mathcal{F}_{1}\left(\bar{t},\hat{G},\tilde{\theta},\Gamma,\dfrac{1}{\omega}\right)\,, \label{eq:dotHatGav_event_4_siso} \\
\frac{d\tilde{\theta}(\bar{t})}{d\bar{t}}&= \frac{1}{\omega}\mathcal{F}_{2}\left(\bar{t},\hat{G},\tilde{\theta},\Gamma,\dfrac{1}{\omega}\right)\,, \label{eq:dotTildeTheta_3_event_siso} \\
\frac{d\Gamma(\bar{t})}{d\bar{t}}&= \frac{1}{\omega}\mathcal{F}_{3}\left(\bar{t},\hat{G},\tilde{\theta},\Gamma,\dfrac{1}{\omega}\right)\,, \label{eq:dotGamma_1_event_siso}
\end{align}
}
with
\begin{align}
&\mathcal{F}_1\left(\bar{t},\hat{G},\tilde{\theta},\dfrac{1}{\omega}\right)\mathbb{=} -\frac{a^{2}}{2}\left(1-\cos\left(2\omega \bar{t}\right)\right)H^{\ast}K\Gamma(\bar{t})\hat{G}(\bar{t})\nonumber \\
&-\frac{a^{2}}{2}\left(1-\cos\left(2\omega \bar{t}\right)\right)H^{\ast}Ke(\bar{t})-\frac{3a^{2}\omega H^{\ast}}{8}\cos\left(3 \omega \bar{t}\right) \nonumber \\
&+a^{2}\omega H^{\ast}\sin\left(2\omega \bar{t}\right)\tilde{\theta}(\bar{t})+\left(a\omega Q^{\ast}+\frac{3a^{2}\omega H^{\ast}}{8}\right)\cos\left(\omega \bar{t}\right)\,, \label{eq:hatMathcalG_siso}
\end{align}
\begin{align}
&\mathcal{F}_2\left(\bar{t},\hat{G},\tilde{\theta},\dfrac{1}{\omega}\right)\mathbb{=}-\frac{a^2 }{2}\left(1-\cos\left(2\omega \bar{t}\right)\right)K\Gamma(\bar{t})H^{\ast}\tilde{\theta}(\bar{t})-Ke(\bar{t})\nonumber \\
&+\left[-\left(a Q^{\ast}+\frac{3 a^{3}H^{\ast}}{8}\right)\sin\left(\omega \bar{t}\right)+\frac{a^{3}H^{\ast}}{8}\sin\left(3\omega \bar{t}\right)\right]K\Gamma(\bar{t})\,, \label{eq:tildeMathcalTheta_siso} \\
&\mathcal{F}_3\left(\bar{t},\hat{G},\tilde{\theta},\dfrac{1}{\omega}\right)\mathbb{=}\omega_{r}\Gamma\left(t\right) - \omega_{r}\hat{H}\left(t\right)\Gamma^2\left(t\right) \,, \quad \omega_{r}>0\,. \label{eq:mathcalGamma_v1} 
\end{align}


Now, by defining the augmented state
\textcolor{black}{
\begin{align} \label{cartagena}
X^{T}(\bar{t}):=\begin{bmatrix}\hat{G}(\bar{t})\,, \tilde{\theta}(\bar{t})\,, \Gamma(\bar{t}) 
\end{bmatrix}\,,
\end{align}
}
one arrives at the dynamics
\textcolor{black}{
\begin{align}
\dfrac{dX(\bar{t})}{d\bar{t}}&=\dfrac{1}{\omega}\mathcal{F}\left(\bar{t},X,\dfrac{1}{\omega}\right)\,, \quad \label{eq:dotX_event_siso}
\mathcal{F}^{T}=\begin{bmatrix} \mathcal{F}_{1}\,, \mathcal{F}_{2}\,, \mathcal{F}_{3} \end{bmatrix}\,.
\end{align}}

Due to the discontinuous nature of the proposed control strategy, the averaging theory for discontinuous systems is employed in the paper, according to \citep{P:1979}.

The augmented system (\ref{eq:dotX_event_siso}) has a small parameter $1/\omega$ as well as a $T$-periodic function $\mathcal{F}\left(\bar{t},X,\dfrac{1}{\omega}\right)$ in $\bar{t}$, hence it can be studied through the averaging method for stability analysis by averaging  $\mathcal{F}\left(\bar{t},X,\dfrac{1}{\omega}\right)$ at $\displaystyle \lim_{\omega\to \infty}\dfrac{1}{\omega}=0$, as shown in \citep{P:1979}, {\it i.e.},
\begin{align}
\dfrac{dX_{\text{av}}(\bar{t})}{d\bar{t}}&=\dfrac{1}{\omega}\mathcal{F}_{\text{av}}\left(X_{\text{av}}\right) \,, \label{eq:dotXav_event_1_siso} \\
\mathcal{F}_{\text{av}}\left(X_{\text{av}}\right)&=\dfrac{1}{T}\int_{0}^{T}\mathcal{F}\left(\delta,X_{\text{av}},0\right)d\delta
\,.  \label{eq:mathcalFav_event_siso}
\end{align}
Basically, the problem in the averaging method is to determine in what sense the behavior of the autonomous system (\ref{eq:dotXav_event_1_siso}) approximates the behavior of the nonautonomous system (\ref{eq:dotX_event_siso}) such that (\ref{eq:dotX_event_siso}) can be represented as a perturbed version of the system (\ref{eq:dotXav_event_1_siso}).

Therefore, treating the non-periodic states $\hat{G}(\bar{t})$, $\tilde{\theta}(\bar{t})$, $\Gamma(\bar{t})$, and $e(\bar{t})$ as constants in (\ref{eq:dotHatGav_event_4_siso})--(\ref{eq:dotX_event_siso}), one gets the following average closed-loop system   
\begin{align}
\frac{d\hat{G}_{\text{av}}(\bar{t})}{d\bar{t}}&=-\frac{a^{2}}{2\omega}H^{\ast}K\Gamma_{\rm{av}}(\bar{t})\hat{G}_{\rm{av}}(\bar{t})-\frac{a^{2}}{2\omega}H^{\ast}Ke_{\rm{av}}(\bar{t})\,, \label{eq:dotHatGav_event_1_siso} \\
\frac{d\tilde{\theta}_{\text{av}}(\bar{t})}{d\bar{t}}&=-\frac{a^2 }{2\omega}K\Gamma_{\rm{av}}(\bar{t})H^{\ast}\tilde{\theta}_{\rm{av}}(\bar{t})-\frac{1}{\omega}Ke_{\rm{av}}(\bar{t})\,, \label{eq:dotTildeThetaAv_event_1_siso} \\
\frac{d\Gamma_{\text{av}}(\bar{t})}{d\bar{t}}&=\frac{\omega_{r}}{\omega}\Gamma_{\rm{av}}(\bar{t})-\frac{\omega_{r}}{\omega}H^{\ast}\Gamma_{\rm{av}}^{2}(\bar{t})\,, \label{eq:dotGammaAv_event_1_siso}
\end{align}
with average update error given by 
\begin{align}
e_{\text{av}}(\bar{t})&=\Gamma_{\rm{av}}(\bar{t_{k}})\hat{G}_{\text{av}}(\bar{t}_{k})-\Gamma_{\rm{av}}(\bar{t})\hat{G}_{\text{av}}(\bar{t})\,. \label{eq:Eav_event_1_siso} 
\end{align}

The equilibrium points of (\ref{eq:dotGammaAv_event_1_siso}), obtained by solving $\dot{\Gamma}_{\rm{av}}\left(t\right) =0$, are $\Gamma^{\ast}_{\rm{av}}=0$ and $\Gamma^{\ast}_{\rm{av}}=\frac{1}{H^{\ast}}$. Thus, we define the Hessian's inverse estimation error
\begin{align} 
\tilde{\Gamma}_{\rm{av}}\left(\bar{t}\right) = \Gamma_{\rm{av}}\left(\bar{t}\right) - \frac{1}{H^{\ast}} \,, \label{eq:GammaTilde_v1} 
\end{align} 
and the system (\ref{eq:dotHatGav_event_1_siso})--(\ref{eq:dotGammaAv_event_1_siso}) can be rewritten with respect to $\tilde{\Gamma}_{\rm{av}}\left(t\right)$ as
\begin{align}
\frac{d\hat{G}_{\text{av}}(\bar{t})}{d\bar{t}}&=-\frac{a^{2}}{2\omega}K\hat{G}_{\rm{av}}(\bar{t})-\frac{a^{2}}{2\omega}H^{\ast}Ke_{\rm{av}}(\bar{t}) \nonumber \\
&\quad-\frac{a^{2}}{2\omega}H^{\ast}K\tilde{\Gamma}_{\rm{av}}(\bar{t})\hat{G}_{\rm{av}}(\bar{t})\,, \label{eq:dotHatGav_20250115_1_siso} \\
\frac{d\tilde{\theta}_{\text{av}}(\bar{t})}{d\bar{t}}&=-\frac{a^2 }{2\omega}K\tilde{\theta}_{\rm{av}}(\bar{t})-\frac{1}{\omega}Ke_{\rm{av}}(\bar{t}) \nonumber \\
&\quad-\frac{a^2 }{2\omega}KH^{\ast}\tilde{\Gamma}_{\rm{av}}(\bar{t})\tilde{\theta}_{\rm{av}}(\bar{t})\,, \label{eq:dotTildeThetaAv_20250115_1_siso} \\
\dot{\tilde{\Gamma}}_{\rm{av}}\left(\bar{t}\right)& = -\frac{\omega_{r}}{\omega}\tilde{\Gamma}_{\rm{av}}\left(\bar{t}\right) - \frac{\omega_{r}}{\omega}H^{\ast}\tilde{\Gamma}^2_{\rm{av}}\left(\bar{t}\right)\,. \label{eq:dotTildeGamma_20250115_1_siso} 
\end{align} 

The system (\ref{eq:dotHatGav_20250115_1_siso})--(\ref{eq:dotTildeGamma_20250115_1_siso}) is nonlinear due the quadratic terms $\tilde{\Gamma}_{\rm{av}}(\bar{t})\hat{G}_{\rm{av}}(\bar{t})$, $\tilde{\Gamma}_{\rm{av}}(\bar{t})\tilde{\theta}_{\rm{av}}(\bar{t})$ and $\tilde{\Gamma}^2_{\rm{av}}\left(\bar{t}\right)$. To analyze its stability properties, we provide a linearization to approximate of the system's behavior in the vicinity of origin. In this sense, for small values of $\tilde{\theta}_{\rm{av}}(0)$ and $\tilde{\Gamma}_{\rm{av}}\left(0\right)$, the quadratic terms $\tilde{\Gamma}_{\rm{av}}(\bar{t})\hat{G}_{\rm{av}}(\bar{t})$, $\tilde{\Gamma}_{\rm{av}}(\bar{t})\tilde{\theta}_{\rm{av}}(\bar{t})$ and $\tilde{\Gamma}^2_{\rm{av}}\left(\bar{t}\right)$ are negligible in linearized corresponding system such that the system (\ref{eq:dotHatGav_20250115_1_siso})--(\ref{eq:dotTildeGamma_20250115_1_siso}) can be rewritten with respect to $\tilde{\Gamma}_{\rm{av}}\left(\bar{t}\right)$ as
\begin{align}
\frac{d\hat{G}_{\text{av}}(\bar{t})}{d\bar{t}}&=-\frac{a^{2}}{2\omega}K\hat{G}_{\rm{av}}(\bar{t})-\frac{a^{2}}{2\omega}H^{\ast}Ke_{\rm{av}}(\bar{t})\,, \label{eq:dotHatGav_20250115_2_siso} \\
\frac{d\tilde{\theta}_{\text{av}}(\bar{t})}{d\bar{t}}&=-\frac{a^2 }{2\omega}K\tilde{\theta}_{\rm{av}}(\bar{t})-\frac{1}{\omega}Ke_{\rm{av}}(\bar{t}) \,, \label{eq:dotTildeThetaAv_20250115_2_siso} \\
\dot{\tilde{\Gamma}}_{\rm{av}}\left(\bar{t}\right)& = -\frac{\omega_{r}}{\omega}\tilde{\Gamma}_{\rm{av}}\left(\bar{t}\right)\,. \label{eq:dotTildeGamma_20250115_2_siso} 
\end{align} 

Hence, from (\ref{eq:dotHatGav_20250115_2_siso})--(\ref{eq:dotTildeGamma_20250115_2_siso}), it is easy to verify the ISS relations of $\hat{G}_{\text{av}}(\bar{t})$ and $\tilde{\theta}_{\text{av}}(\bar{t})$ with respect to the averaged measurement error  $e_{\text{av}}(\bar{t})$. 

Moreover, from (\ref{eq:hatG_4_event_siso}), 
\begin{align}
\hat{G}_{\text{av}}(\bar{t})= \frac{a^{2}H^{\ast}}{2}\tilde{\theta}_{\text{av}}(\bar{t})\,, \label{eq:hatGav_event_1_siso}
\end{align}
and, consequently,
\begin{align}
\tilde{\theta}_{\text{av}}(\bar{t})= \frac{2}{a^{2}H^{\ast}}\hat{G}_{\text{av}}(\bar{t})\,. \label{eq:tildeThetaAv_event_1_siso}
\end{align}

Therefore, the following average event-triggered equivalences can be introduced for the average system.

\begin{definition} \label{def:averageStaticEvent} Consider the nonlinear mapping $\Xi: \mathbb{R} \times \mathbb{R}\mapsto \mathbb{R}$ given by 
\begin{align}
\Xi(\hat{G}_{\rm{av}},e_{\rm{av}}) = \sigma|\hat{G}_{\rm{av}}(\bar{t})|-\beta |e_{\rm{av}}(\bar{t})|\,,  \label{eq:Psi_event_1_siso}  
\end{align}
with $K>0$ being the control gain in (\ref{eq:U_event_siso}). The event-triggered controller with triggering condition consists of the next two components:
\begin{enumerate}
	\item \textcolor{black}{A sequence of increasing times} $I=\{\bar{t}_{0}\,, \bar{t}_{1}\,, \bar{t}_{2}\,,\ldots\}$, with $\bar{t}_{0}=0$, generated under the following rules 
		\begin{itemize}
			\item If $\left\{\bar{t} \in\mathbb{R}^{+}: \bar{t}>\bar{t}_{k} ~ \wedge ~ \Xi (\hat{G}_{\rm{av}},e_{\rm{av}}) > 0 \right\} = \emptyset$, then the set of the times of the events is $I=\{\bar{t}_{0}\,, \bar{t}_{1}\,, \ldots, \bar{t}_{k}\}$.
			\item If $\left\{\bar{t} \in\mathbb{R}^{+}: \bar{t}>\bar{t}_{k} ~ \wedge ~ \Xi (\hat{G}_{\rm{av}},e_{\rm{av}}) < 0 \right\} \neq \emptyset$, the next event time is given by
				\begin{align}
					\bar{t}_{k+1}&=\textcolor{black}{\inf}\left\{\bar{t} \in\mathbb{R}^{+}: \bar{t}\!>\!\bar{t}_{k} ~ \!\wedge\! ~ \Xi (\hat{G}_{\rm{av}},e_{\rm{av}}) \!<\!0 \right\}\,, \label{eq:tk+1_event_av_siso}
				\end{align}
				which is the average event-trigger mechanism.
		\end{itemize}
	\item A feedback control action updated at the generated triggering instants given by
		\begin{align}
			u^{\rm{av}}_{k}\!=\!-K\Gamma_{\rm{av}}(\bar{t}_{k})\hat{G}_{\rm{av}}(\bar{t}_{k}) \,, \quad \forall \bar{t} \!\in\! \lbrack \bar{t}_{k}\,, \bar{t}_{k+1}\phantom{(}\!\!)\,, \quad k\!\in\!\mathbb{N}\,. \label{eq:U_MD3}
		\end{align}
\end{enumerate}  
\end{definition}

\section{Stability Analysis} \label{sec:stblt_siso}

This section does not assumes any knowledge of the nonlinear map (\ref{eq:y_event_siso}). 
%
The following assumptions are considered in the main theorem:

\begin{description}
\item[{\bf(A1)}] The unique optimizer value $\theta^{\ast} \in R$ and the scalar extremum point $Q^{\ast}$ are unknown parameters of the nonlinear map (\ref{eq:y_event_siso}).
\item[{\bf(A2)}] The Hessian $H^{\ast}$ is an unknown parameter.
\item[{\bf(A3)}] The control gain satisfies $K>0$.
%
\item[{\bf(A4)}] There is no delay due processing of sensor and actuator as well as transmission in the sensor-to-controller and controller-to-actuator branches.
\item[{\bf(A5)}] Only $\hat{G}(t)$ is available for the event-triggered design.     
\end{description}
\medskip

\begin{theorem} \label{thm:NETESC_1_siso}
Consider the closed-loop average system (\ref{eq:dotHatGav_20250115_2_siso})--(\ref{eq:dotTildeGamma_20250115_2_siso}) and the average SET-NewtonES mechanism given by (\ref{eq:tk+1_event_av_siso}). Suppose that Assumptions (A1)--(A5) are hold. If $\Xi(\hat{G}_{\rm{av}},e_{\rm{av}})$ is given by (\ref{eq:Psi_event_1_siso}) and $\omega>0$ in (\ref{eq:omega_event_1_siso}) is a constant sufficiently large, then, \textcolor{black}{for $\tilde{\theta}(0)$ and $\tilde{\Gamma}(0)$  sufficiently small}, 
the average system (\ref{eq:dotHatGav_20250115_2_siso})--(\ref{eq:dotTildeGamma_20250115_2_siso}) is locally exponentially stable and the following relations can be obtained for the closed-loop system (\ref{eq:dotX_event_siso}), with state in (\ref{cartagena}), in the original variables $\tilde{\theta}(t)\,, \Gamma(t)$:
\begin{small}
\begin{align}
|\theta(t)-\theta^{\ast}|&\leq \exp\left(\mathbb{-}\frac{(1\mathbb{-}\sigma) a^{2}K}{2}t\right) |\theta(0)\mathbb{-}\theta^{\ast}|\mathbb{+}\mathcal{O}\left(a\mathbb{+}\frac{1}{\omega}\right), \label{eq:normTheta_thm1_siso} \\ 
|y(t)\mathbb{-}Q^{\ast}|&\leq \exp\left(\!\mathbb{-}\frac{(1\mathbb{-}\sigma) a^{2}K}{2}t\right) |y(0)\mathbb{-}Q^{\ast}|\mathbb{+}\mathcal{O}\left(a^2\mathbb{+}\frac{1}{\omega^2}\right), \label{eq:normY_thm1_siso}\\
\left|\Gamma(t)\mathbb{-}\frac{1}{H^{\ast}}\right|&\leq \exp\left(\mathbb{-}\omega_{r}t\right) \left|\Gamma(0)\mathbb{-}\frac{1}{H^{\ast}}\right|\mathbb{+}\mathcal{O}\left(\frac{1}{\omega}\right), \quad a>0 \,, \ \omega_r>0. \label{eq:normGamma_thm1_siso}
\end{align}
\end{small}
\end{theorem}

\begin{pf} 
See Appendix~\ref{goldenglobe}. \hfill $\square$
\end{pf}

%

\section{Simulation results} \label{sec:sim_siso}

To highlight the main ideas of the proposed SET-NewtonES control strategy, the nonlinear map (\ref{eq:y_event_siso}) has the following unknown parameters: 
$H^{\ast}=-0.15$, $Q^{\ast}=7$ and $\theta^{\ast}=5$. The parameters of the probing and demodulation signals are:  $a=0.1$, $\omega=3$ [rad/sec], and we select the event-triggered parameter $\sigma=0.9$. The control gain is $K=18$, and the initial conditions are $\hat{\theta}(0)=2.0$, $\Gamma(0)=-0.1$. For the Riccati filter parameter, we have  
$\omega_r$ = 1.0 [rad/s]. The following plots consider a numerical simulation over a period from $0$ to $500$ seconds.

\begin{figure*}[h!]
	\centering
	\subfigure[Nonlinear map input, $\theta(t)$. \label{fig:theta}]{\includegraphics[width=6.3cm]{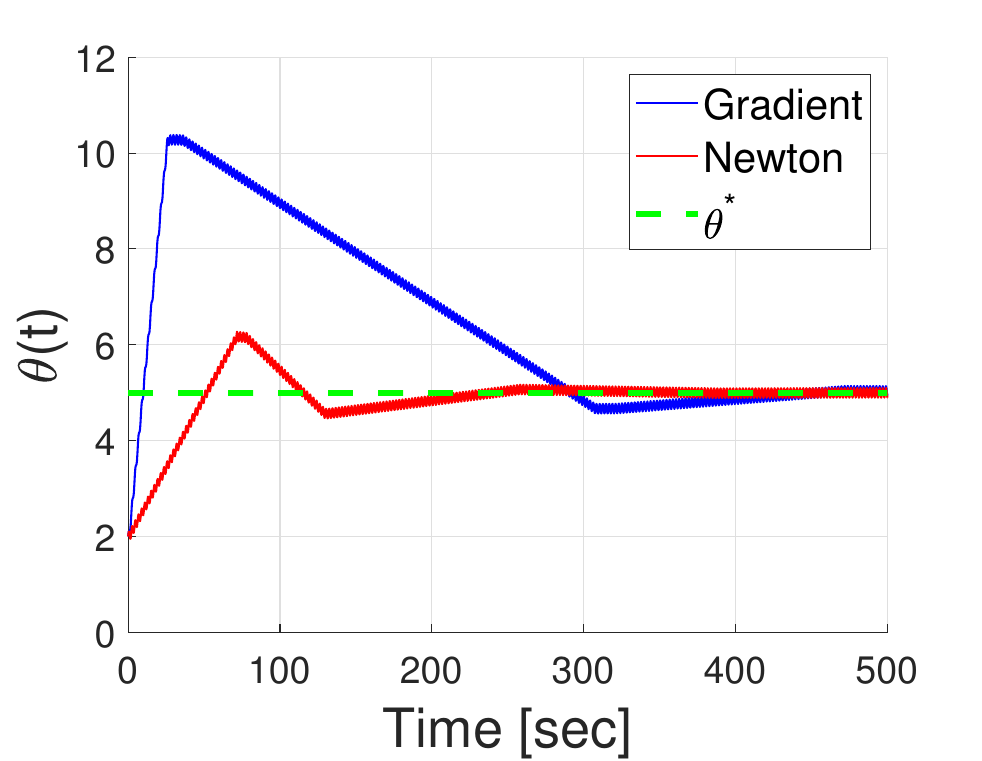}}
	\subfigure[Control input, $U(t)$. \label{fig:U}]{\includegraphics[width=6.3cm]{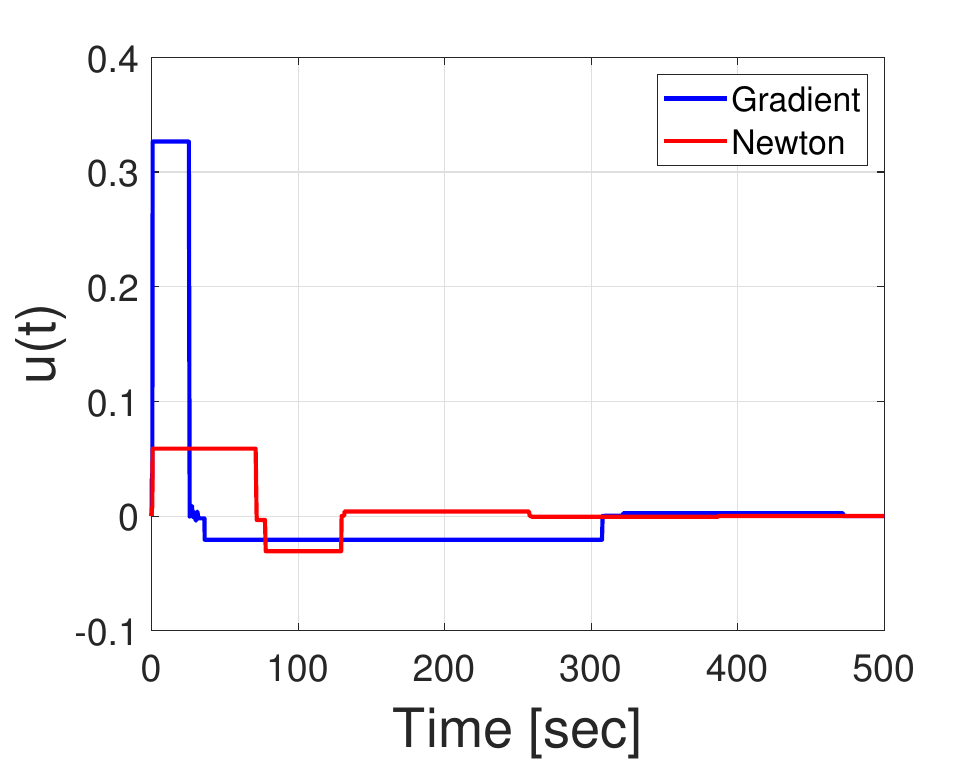}}
		\\
		\subfigure[Controller update. \label{fig:update}]{\includegraphics[width=6.3cm]{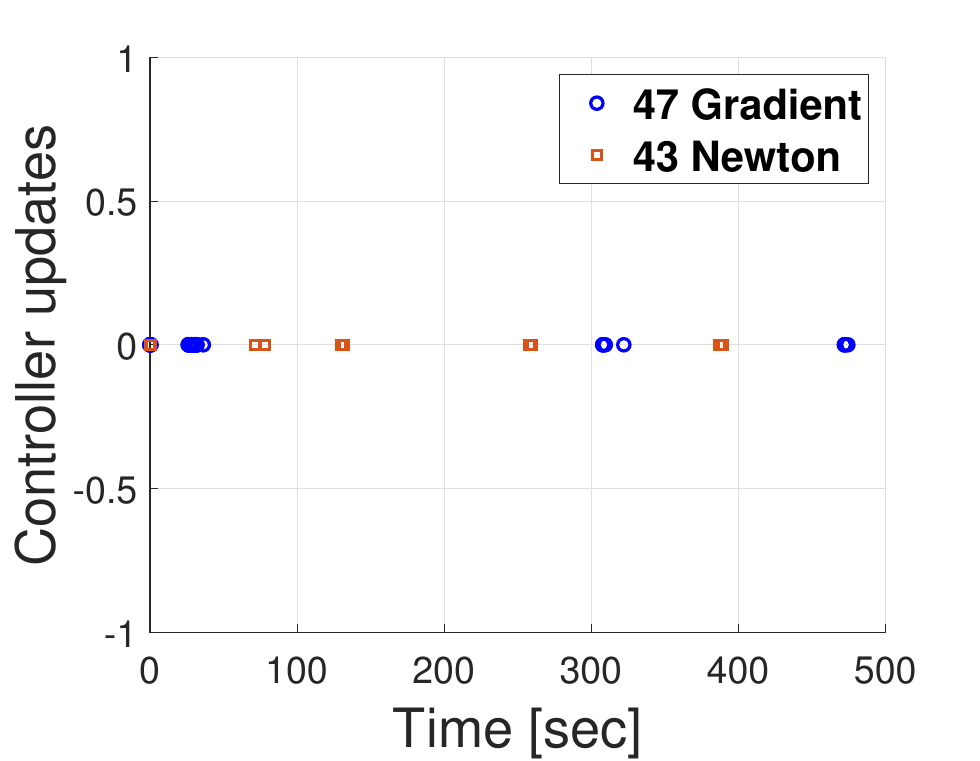}}
		\subfigure[Hessian inverse $\Gamma(t)$. \label{fig:gamma}]{\includegraphics[width=6.3cm]{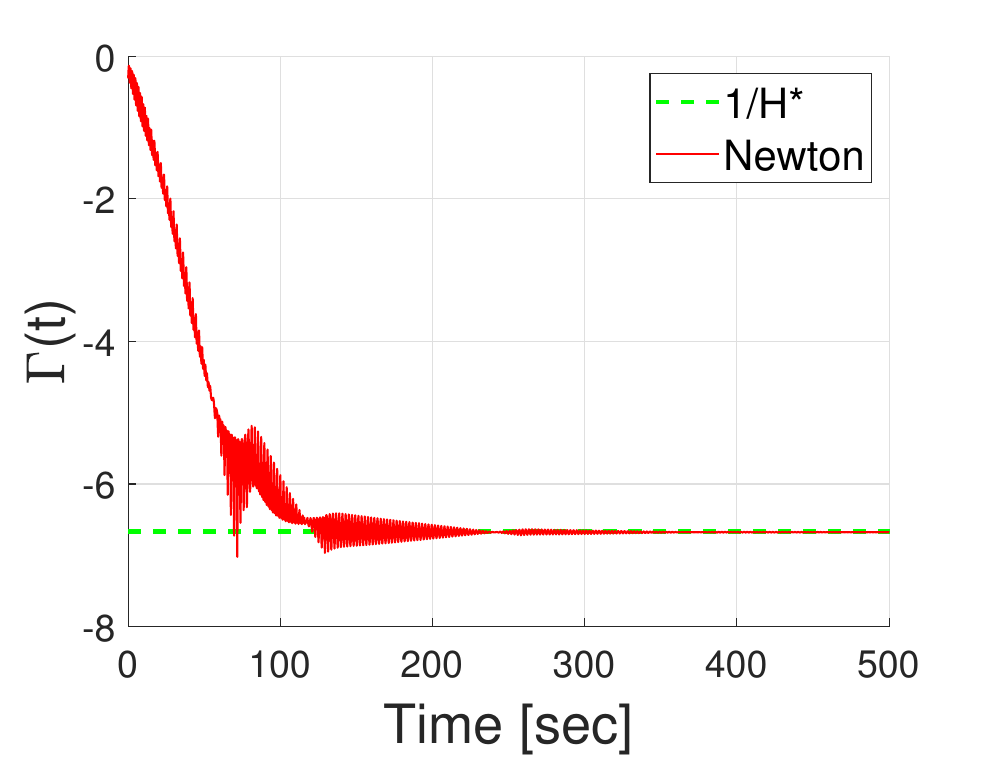}}
\end{figure*}

 Fig.~\ref{fig:theta} illustrates the convergence of the input $\theta(t)$ to the optimizer $\theta^*$ of the static map, the blue curve represents the response of the SET-GradientES  method of Fig.~\ref{fig:BD_SET_GradientES_SISO}, the red curve represents the response of the proposed SET-NewtonES method, and the green curve marks the desired optimum, $\theta^*=7$. The variable $\theta(t)$ from Fig.~\ref{fig:theta} oscillates around the optimal value after $100$ seconds by using the Newton-based method, while the gradient-based method still remains in its transient phase. Fig.~\ref{fig:gamma}   
 for $\Gamma(t)$ 
 highlights that a faster convergence rate with the Newton-based method compared to the gradient-based method and can only be achieved due to the estimate of the Hessian's inverse $\Gamma^* = 1/H^{*}=-6.7$.  Furthermore, the gradient-based method exhibits higher amplitudes in its trajectory, which is not seen in the Newton-based control signal, see Fig. \ref{fig:U}. Fig.~\ref{fig:update} displays the control signal updating times, showing that the gradient-based method updated the controller $47$ times, whereas the Newton-based method required only $43$ updates. This demonstrates that the Newton-based approach not only achieves the optimal value with fewer control signal updates but also improves the convergence rate.

\section{Conclusion} \label{sec:concl_siso}

In this paper, we have proposed a comparison between the Newton-based extremum seeking and the gradient-based method, both employing an event-triggered mechanism in order to optimize scalar-static nonlinear maps. It is shown that the convergence rate of the proposed Newton-based extremum seeking is faster than the companion gradient scheme, even when the event-triggered methodology is incorporated in the feedback loop. By using a Riccati equation filter to compute the inverse of the unknown Hessian of the map, the convergence rate can be arbitrarily chosen by the user. Additionally, the feedback update action occurs less frequently than in the gradient method, thereby reducing computational resource usage.



\appendix
\section{Proof of Theorem~\ref{thm:NETESC_1_siso}} \label{goldenglobe}   
\begin{small}
The proof of the theorem is divided into two parts: stability analysis and avoidance of Zeno behavior.

\begin{flushleft}
\underline{A. Stability Analysis}
\end{flushleft}
\medskip

Consider the following Lyapunov function candidate
\begin{align}
V_{\text{av}}(\bar{t})=\hat{G}^{2}_{\text{av}}(\bar{t}) \,, \label{eq:Vav_event_siso}
\end{align}  
with time-derivative, by using (\ref{eq:dotHatGav_20250115_2_siso}), given by
\begin{align}
\frac{dV_{\text{av}}(\bar{t})}{d\bar{t}}&=-\frac{a^{2}}{\omega}K\hat{G}_{\rm{av}}^{2}(\bar{t})-\frac{a^{2}}{\omega}H^{\ast}K\hat{G}_{\rm{av}}(\bar{t})e_{\rm{av}}(\bar{t}) \nonumber \\ 
 &\leq-\frac{a^{2}}{\omega}K\hat{G}_{\rm{av}}^{2}(\bar{t})+\frac{a^{2}}{\omega}\beta K|\hat{G}_{\rm{av}}(\bar{t})||e_{\rm{av}}(\bar{t})|\,, \beta>|H^{\ast}|\,, \nonumber \\
 &\leq-\frac{a^{2}K}{\omega}|\hat{G}_{\rm{av}}(\bar{t})|\left(|\hat{G}_{\rm{av}}(\bar{t})|-\beta |e_{\rm{av}}(\bar{t})|\right) \,. \label{eq:dotVav_event_1_siso}
\end{align} 
From (\ref{eq:dotVav_event_1_siso}), if there is no measurement error $e(t)$, {\it i.e.}, $e(t)\equiv 0$ $\forall t>0$, therefore the classic extremum seeking implementation, according to equation (\ref{eq:dotVav_event_1_siso}) becomes $\frac{dV_{\text{av}}}{d\bar{t}}= -\frac{a^{2}K}{\omega}\hat{G}^{2}_{\text{av}}(\bar{t})$. On the other hand, the proposed event-triggered approach, with update law (\ref{eq:tk+1_event_av_siso}) and $\Xi(\hat{G}_{\rm{av}},e_{\rm{av}})$ given by (\ref{eq:Psi_event_1_siso}), ensures that, between to events, $|e_{\rm{av}}(\bar{t})|\leq \dfrac{\sigma}{\beta}|\hat{G}_{\rm{av}}(\bar{t})|$ and, consequently, in closed-loop an exponential decay of $V_{\text{av}}(\bar{t})$ given by a pre-specified fraction of the ideal decay rate such that
\begin{align}
\frac{dV_{\text{av}}(\bar{t})}{d\bar{t}}\leq-\frac{(1-\sigma) a^{2}K}{\omega}\hat{G}^{2}_{\text{av}}(\bar{t})\,. \label{eq:dotVav_event_2_siso}
\end{align}

The event-triggered mechanism supervises the time derivative of the Lyapunov function given by (\ref{eq:dotVav_event_1_siso}) and its pre-specified upper bound in (\ref{eq:dotVav_event_2_siso}) to set the instant on which these signals meet. This time-instant is the same to send data over the network and update the actuator, where the condition $\Xi(\hat{G}_{\text{av}},e_{\text{av}})<0$ is verified as in (\ref{eq:tk+1_event_av_siso}). This process can \textcolor{black}{take} place for an indefinite number of times, in other words, whenever necessary, and guarantees the exponential stability of $\hat{G}_{\text{av}}(\bar{t})$ in closed-loop.

By using (\ref{eq:Vav_event_siso}) and (\ref{eq:dotVav_event_2_siso}), for $\bar{t} \in [\!\phantom{]}\bar{t}_{k}\,,\bar{t}_{k+1})$, an upper bound for (\ref{eq:dotVav_event_1_siso}) is 
\begin{align}
\frac{dV_{\text{av}}(\bar{t})}{d\bar{t}}&\leq -\frac{(1-\sigma) a^{2}K}{\omega}V_{\text{av}}(\bar{t}) \,. \label{eq:dotVav_event_3_siso}
\end{align}

Then, invoking the Comparison Lemma \citep{K:2002} an upper bound $\bar{V}_{\text{av}}(\bar{t})$ for $V_{\text{av}}(\bar{t})$ is
\begin{align}
V_{\text{av}}(\bar{t})\leq \bar{V}_{\text{av}}(\bar{t}) \,, \quad \forall \bar{t}\in \lbrack \bar{t}_{k},\bar{t}_{k+1}\phantom{(}\!\!)\,. \label{eq:VavBarVav_1_siso}
\end{align}
given by the solution of the following dynamics
\begin{align}
\frac{d\bar{V}_{\text{av}}(\bar{t})}{d\bar{t}}=-\frac{(1-\sigma) a^{2}K}{\omega}\bar{V}_{\text{av}}(\bar{t}),~~\bar{V}_{\text{av}}(\bar{t}_{k})=V_{\text{av}}(\bar{t}_{k})\,,
\end{align}
In other words, $ \forall \bar{t}\in \lbrack \bar{t}_{k},\bar{t}_{k+1}\phantom{(}\!\!)$,
\begin{align}
\bar{V}_{\text{av}}(\bar{t})=\exp\left(-\frac{(1-\sigma) a^{2}K}{\omega}(\bar{t}-\bar{t}_{k})\right)V_{\text{av}}(\bar{t}_{k})\,, \label{eq:_siso}
\end{align}
and the inequality (\ref{eq:VavBarVav_1_siso}) is rewritten as
\begin{align}
V_{\text{av}}(\bar{t})\leq \exp\left(-\frac{(1-\sigma) a^{2}K}{\omega}(\bar{t}-\bar{t}_{k})\right)V_{\text{av}}(\bar{t}_{k}) \,. \label{eq:VavBarVav_2_siso}
\end{align}
By defining, $\bar{t}_{k}^{+}$ and $\bar{t}_{k}^{-}$ as the right and left limits of $\bar{t}=\bar{t}_{k}$, respectively, it easy to verify that $$V_{\text{av}}(\bar{t}_{k+1}^{-})\leq \exp\left(-\dfrac{(1-\sigma)a^{2}K}{\omega}(\bar{t}_{k+1}^{-}-\bar{t}_{k}^{+})\right)V_{\text{av}}(\bar{t}_{k}^{+}).$$ Since $V_{\text{av}}(\bar{t})$ is continuous, $V_{\text{av}}(\bar{t}_{k+1}^{-})=V_{\text{av}}(\bar{t}_{k+1})$ and $V_{\text{av}}(\bar{t}_{k}^{+})=V_{\text{av}}(\bar{t}_{k})$, and therefore,
\begin{align}
    V_{\text{av}}(\bar{t}_{k+1})\leq \exp\left(-\frac{(1-\sigma)a^{2}K}{\omega}(\bar{t}_{k+1}-\bar{t}_{k})\right)V_{\text{av}}(\bar{t}_{k})\,. \label{eq:mmd_1_s}
\end{align}
Hence, for any $\bar{t}\geq 0$ in $ \bar{t}\in \lbrack \bar{t}_{k},\bar{t}_{k+1}\phantom{(}\!\!)$, $k \in \mathbb{N}$, one has 
\begin{align}
    &V_{\text{av}}(\bar{t})\leq \exp\left(-\frac{(1-\sigma) a^{2}K}{\omega}(\bar{t}-\bar{t}_{k})\right) V_{\text{av}}(\bar{t}_{k}) \nonumber \\
    &\leq \exp\left(-\frac{(1-\sigma) a^{2}K}{\omega}(\bar{t}-\bar{t}_{k})\right) \times \nonumber \\
    &\quad \times \exp\left(-\frac{(1-\sigma) a^{2}K}{\omega}(\bar{t}_{k}-\bar{t}_{k-1})\right)V_{\text{av}}(\bar{t}_{k-1}) \nonumber \\
    &\leq \ldots \leq \nonumber \\
    &\leq \exp\left(-\frac{(1-\sigma) a^{2}K}{\omega}(\bar{t}-\bar{t}_{k})\right)\times \nonumber \\
    &\times \prod_{i=1}^{i=k}\exp\left(-\frac{(1-\sigma) a^{2}K}{\omega}(\bar{t}_{i}-\bar{t}_{i-1})\right)V_{\text{av}}(\bar{t}_{i-1}) \nonumber \\
    &=\exp\left(-\frac{(1-\sigma) a^{2}K}{\omega}\bar{t}\right) V_{\text{av}}(0)\,.
\end{align}

Now, one obtains
\begin{align}
|\hat{G}_{\text{av}}(\bar{t})|^{2}&\leq \exp\left(-\frac{(1-\sigma) a^{2}K}{\omega}\bar{t}\right) |\hat{G}_{\text{av}}(0)|^{2} \label{eq:VavBarVav_3_siso} \\
&=\left[\exp\left(-\frac{(1-\sigma)a^{2}K}{2\omega}\bar{t}\right)|\hat{G}_{\text{av}}(0)|\right]^{2}\,. \label{eq:VavBarVav_4_siso}
\end{align}
Hence,
\begin{align}
|\hat{G}_{\text{av}}(\bar{t})|\leq \exp\left(-\frac{(1-\sigma) a^{2}K}{2\omega}\bar{t}\right)|\hat{G}_{\text{av}}(0)|\,. \label{eq:normHatGav_1_siso}
\end{align}
Although the analysis has been focused on the convergence of $\hat{G}_{\text{av}}(\bar{t})$ and, consequently $\hat{G}(t)$, the obtained results through (\ref{eq:hatG_6_siso}) can be easily extended to the variable $\tilde{\theta}_{\text{av}}(\bar{t})$ and $\tilde{\theta}(t)$. Notice that, by using (\ref{eq:hatGav_event_1_siso}), inequality (\ref{eq:normHatGav_1_siso}) can be rewritten in terms of $\tilde{\theta}_{\text{av}}(\bar{t})$, {\it i.e.},
\begin{align}
|\tilde{\theta}_{\text{av}}(\bar{t})|\leq\exp\left(-\frac{(1-\sigma) a^{2}K}{2\omega}\bar{t}\right)|\tilde{\theta}_{\text{av}}(0)|\,. \label{eq:normTildeThetaAv_1_siso}
\end{align}
Moreover, the solution of (\ref{eq:dotTildeGamma_20250115_2_siso}) satisfies 
\begin{align}
|\tilde{\Gamma}_{\rm{av}}(\bar{t})|\leq \exp\left(-\frac{\omega_{r}}{\omega}\bar{t}\right)|\tilde{\Gamma}_{\rm{av}}(0)|\,. \label{eq:normGammaAv_1_siso}
\end{align}
Since (\ref{eq:Gamma}),(\ref{eq:dotHatGav_event_3_siso}) and (\ref{eq:dotTildeTheta_2_event_siso}) are $T$-periodic in $t$, $1/\omega$ is a positive small parameter, and from inequalities (\ref{eq:normHatGav_1_siso}) and (\ref{eq:normTildeThetaAv_1_siso}) the origin $\hat{G}_{\text{av}}=\tilde{\theta}_{\text{av}}=\tilde{\Gamma}_{\rm{av}}=0$ is at least an exponentially stable equilibrium point of the closed-loop event-triggered system. Then, by invoking  \citep[Theorem~2]{P:1979}, there exists an upper bound for (\ref{eq:hatG_4_event_siso}) such that
\begin{align}
|\hat{G}(\bar{t})|&\mathbb{\leq}|\hat{G}_{\text{av}}(\bar{t})|\mathbb{+}\mathcal{O}\left(\frac{1}{\omega}\right)
\mathbb{\leq} \exp\left(\!\!\mathbb{-}\frac{(1\mathbb{-}\sigma) a^{2}K}{2\omega}\bar{t}\right)|\hat{G}_{\text{av}}(0)|\mathbb{+}\mathcal{O}\left(\frac{1}{\omega}\right)\,, \label{eq:hatG_6_siso} \\
|\tilde{\theta}(\bar{t})|&\mathbb{\leq}|\tilde{\theta}_{\text{av}}(\bar{t})|\mathbb{+}\mathcal{O}\left(\frac{1}{\omega}\right)
\mathbb{\leq} \exp\left(\!\!\mathbb{-}\frac{(1\mathbb{-}\sigma) a^{2}K}{2\omega}\bar{t}\right)|\tilde{\theta}_{\text{av}}(0)|\mathbb{+}\mathcal{O}\left(\frac{1}{\omega}\right)\,, \label{eq:nomrTildeTheta_siso} \\
|\tilde{\Gamma}(\bar{t})|&\leq|\tilde{\Gamma}_{\text{av}}(\bar{t})|+\mathcal{O}\left(\frac{1}{\omega}\right)\leq \exp\left(-\frac{\omega_{r}}{\omega}\bar{t}\right)|\tilde{\Gamma}_{\rm{av}}(0)|+\mathcal{O}\left(\frac{1}{\omega}\right)\,. \label{eq:normTildeGamma_siso}
\end{align}

Now, from (\ref{eq:theta_2_event_siso}), we have
\begin{align}
\theta(t)-\theta^{\ast}=\tilde{\theta}(t)+a \sin(\omega t)\,, \label{eq:theta_3_event_siso}
\end{align}
whose norm satisfies 
\begin{align}
|\theta(t)-\theta^{\ast}|&=|\tilde{\theta}(t)+a \sin(\omega t)| \leq |\tilde{\theta}(t)|+|a \sin(\omega t)| \nonumber \\
&\leq \exp\left(\mathbb{-}\frac{(1\mathbb{-}\sigma) a^{2}K}{2}t\right) |\theta(0)\mathbb{-}\theta^{\ast}|\mathbb{+}\mathcal{O}\left(a\mathbb{+}\frac{1}{\omega}\right). \label{eq:theta_4_event_siso}
\end{align}
Defining the error variable $\tilde{y}(t)$ as 
\begin{align}
\tilde{y}(t):=y(t)-Q^{\ast}\,,
\end{align}
and using (\ref{eq:y_event_siso}) as well as the Cauchy-Schwarz Inequality, its absolute value satisfies
\begin{align}
|\tilde{y}(t)|&=|y(t)-Q^{\ast}| =\frac{|H^{\ast}|}{2}|\theta(t)-\theta^{\ast}|^{2}\,, \label{tildeY_event_1_siso}
\end{align}
and its upper bounded with (\ref{eq:theta_4_event_siso})is given by
\begin{align}
|y(t)\mathbb{-}Q^{\ast}| &\leq  \exp\left(\!\!\mathbb{-}\frac{(1\mathbb{-}\sigma) a^{2}K}{2}t\right) |y(0)\mathbb{-}Q^{\ast}|\mathbb{+}\mathcal{O}\left(\! a^2\mathbb{+}\frac{1}{\omega^2} \!\right) \!\!. \label{tildeY_event_4_siso}
\end{align}

Therefore, the inequalities (\ref{eq:normTheta_thm1_siso})--(\ref{eq:normGamma_thm1_siso}) are satisfied.

\begin{flushleft}
\underline{B. Avoidance of Zeno Behavior}
\end{flushleft}
\medskip

Notice that, from (\ref{eq:Psi_event_1_siso}) and (\ref{eq:tk+1_event_av_siso}), and using the Peter-Paul inequality \cite{W:1971}, we can write $cd\leq \frac{c^2}{2\epsilon}+\frac{\epsilon d^2}{2}$, for all $c,d,\epsilon>0$, with $c=|e_{\rm{av}}(\bar{t})|$, $d=|\hat{G}_{\rm{av}}(\bar{t})|$, $\epsilon=\frac{\sigma}{\beta}$ and $\bar{t}\in \lbrack t_{k}\,,t_{k+1}\phantom{(}\!\!\!)$. The following holds
\begin{align}
&\sigma |\hat{G}_{\text{av}}(\bar{t})|^2\mathbb{-}\beta |e_{\text{av}}(\bar{t})||\hat{G}_{\text{av}}(\bar{t})| \geq \nonumber \\
& \sigma|\hat{G}_{\rm{av}}(\bar{t})|^{2}\mathbb{-}\frac{\beta}{2}\!\left(\!\frac{\sigma}{\beta}|\hat{G}_{\rm{av}}(\bar{t})|^2\mathbb{+}\frac{\beta}{\sigma}|e_{\rm{av}}(\bar{t})|^2\!\right)\nonumber \\
&= q|\hat{G}_{\rm{av}}(\bar{t})|^{2}-p|e_{\rm{av}}(\bar{t})|^2\,,\label{ineq:interEvents_2_dynamic_pf4}
\end{align}
where 
\begin{align}
q&=\frac{\sigma}{2} \quad \mbox{and} \quad p=\frac{\beta^2}{2\sigma}\,.\label{ineq:interEvents_3_dynamic_pf4}
\end{align} 
The minimum dwell-time of the event-triggered framework is given by the time it takes for the function
\begin{align}
\phi_{\rm{av}}(\bar{t})=\sqrt{\frac{p}{q}}\frac{|e_{\rm{av}}(\bar{t})|}{|\hat{G}_{\rm{av}}(\bar{t})|} \,, \label{eq:phi_1_dynamic_pf4}
\end{align}
to go from 0 to 1. The derivative of $\phi_{\text{av}}(\bar{t})$ in (\ref{eq:phi_1_dynamic_pf4}) is given by
\begin{align}
\frac{d\phi_{\rm{av}}(\bar{t})}{d\bar{t}}&=\sqrt{\frac{p}{q}}\frac{1}{|e_{\rm{av}}(\bar{t})||\hat{G}_{\rm{av}}(\bar{t})|}\left[e_{\rm{av}}(\bar{t})\frac{de_{\rm{av}}(\bar{t})}{d\bar{t}}\right.\nonumber \\
&\quad \left.-\hat{G}_{\rm{av}}(\bar{t})\frac{d\hat{G}_{\rm{av}}(\bar{t})}{d\bar{t}}\left(\frac{|e_{\rm{av}}(\bar{t})|}{|\hat{G}_{\rm{av}}(\bar{t})|}\right)^2\right]\,. \label{eq:dotPhi_1_static_pf2}
\end{align}

Now, by using (\ref{eq:GammaTilde_v1}), (\ref{eq:dotHatGav_20250115_1_siso}) and (\ref{eq:dotTildeGamma_20250115_1_siso}), the time derivative of (\ref{eq:Eav_event_1_siso}) is 
\begin{align}
&\frac{d e_{\rm{av}}(\bar{t})}{d\bar{t}}=\frac{a^{2}K}{2\omega H^{\ast}}\hat{G}_{\rm{av}}(\bar{t})+\frac{a^{2}K}{2\omega }e_{\rm{av}}(\bar{t})+\frac{a^{2}H^{\ast}K}{2\omega H^{\ast}}\tilde{\Gamma}_{\rm{av}}(\bar{t})\hat{G}_{\rm{av}}(\bar{t}) \nonumber \\
&+\frac{\omega_{r}}{\omega}\tilde{\Gamma}_{\rm{av}}\left(\bar{t}\right)\hat{G}_{\rm{av}}\left(\bar{t}\right)+ \frac{\omega_{r}}{\omega}H^{\ast}\tilde{\Gamma}^2_{\rm{av}}\left(\bar{t}\right)\hat{G}_{\rm{av}}\left(\bar{t}\right)+\frac{a^{2}}{2\omega}K\tilde{\Gamma}_{\rm{av}}\left(t\right)\hat{G}_{\rm{av}}(\bar{t}) \nonumber \\
&+\frac{a^{2}}{2\omega}H^{\ast}K\tilde{\Gamma}_{\rm{av}}\left(\bar{t}\right)e_{\rm{av}}(\bar{t})+\frac{a^{2}}{2\omega}H^{\ast}K\tilde{\Gamma}_{\rm{av}}^{2}(\bar{t})\hat{G}_{\rm{av}}(\bar{t})\,. \label{eq:dotEav_1}
\end{align}

The dynamics  (\ref{eq:dotEav_1}) exhibits a nonlinear behavior due the quadratic and cubic terms between the variables $\tilde{\Gamma}_{\rm{av}}(\bar{t})$ $\hat{G}_{\rm{av}}(\bar{t})$, and $e_{\rm{av}}(\bar{t})$. Thus, we provide a linearization to approximate of the behavior (\ref{eq:dotEav_1}) in the vicinity of origin. This linearization is given by
\begin{align}
\frac{d e_{\rm{av}}(\bar{t})}{d\bar{t}}&=\frac{a^{2}K}{2\omega H^{\ast}}\hat{G}_{\rm{av}}(\bar{t})+\frac{a^{2}K}{2\omega }e_{\rm{av}}(\bar{t})\,. \label{eq:dotEav_2}
\end{align}

Now, using (\ref{eq:dotHatGav_20250115_2_siso}) and (\ref{eq:dotEav_2}),  an upper bound for (\ref{eq:dotPhi_1_static_pf2}) is 
\begin{align}
&\frac{d\phi_{\rm{av}}(\bar{t})}{d\bar{t}}\leq\frac{a^{2}K}{2 \omega}\sqrt{\frac{p}{q}} \left(\frac{1}{|H^{\ast}|}+2\frac{|e_{\rm{av}}(\bar{t})|}{|\hat{G}_{\rm{av}}(\bar{t})|}+|H^{\ast}|\frac{|e_{\rm{av}}(\bar{t})|^2}{|\hat{G}_{\rm{av}}(\bar{t})|^2}\right) \nonumber \\
&\mathbb{\leq}\frac{a^{2}K}{2 \omega}\max\left\{\frac{1}{|H^{\ast}|},1,|H^{\ast}|\right\}\sqrt{\frac{p}{q}}\left(1\mathbb{+}\frac{|e_{\rm{av}}(\bar{t})|}{|\hat{G}_{\rm{av}}(\bar{t})|}\right)^2\,. \label{eq:dotPhi_1_1_static_pf2}
\end{align}

Hence,  using (\ref{eq:phi_1_dynamic_pf4}), inequality (\ref{eq:dotPhi_1_1_static_pf2}) is rewritten as
\begin{align}
&\frac{d\phi_{\rm{av}}(\bar{t})}{d\bar{t}}\leq\frac{a^{2}K}{2 \omega}\max\left\{\frac{1}{|H^{\ast}|},1,|H^{\ast}|\right\}\sqrt{\frac{q}{p}}\left(\sqrt{\frac{p}{q}}+\phi_{\text{av}}(\bar{t})\right)^{2}
\end{align}
Then, by invoking \cite[Comparison Lemma]{K:2002}, an upper bound $\tilde{\phi}_{\text{av}}(\bar{t})$ for $\phi_{\text{av}}(\bar{t})$ according to 
\begin{align}
\phi_{\text{av}}(\bar{t})\leq \tilde{\phi}_{\text{av}}(\bar{t}) \,, \quad \phi_{\text{av}}(0)= \tilde{\phi}_{\text{av}}(0)=0 \,, \quad  \forall \bar{t}\in \lbrack \bar{t}_{k},\bar{t}_{k+1}\phantom{(}\!\!) \,, \label{eq:tildePhi_v2}
\end{align}
is given by the solution of the equation
\begin{align}
\frac{d\tilde{\phi}\textcolor{black}{_{\text{av}}}(\bar{t})}{d\bar{t}}&=\frac{a^{2}K}{2 \omega}\max\left\{\frac{1}{|H^{\ast}|},1,|H^{\ast}|\right\}\sqrt{\frac{q}{p}}\left(\sqrt{\frac{p}{q}}+\tilde{\phi}_{\text{av}}(\bar{t})\right)^{2}\,. \label{eq:dotTildePhi_v2}
\end{align}
The solution of (\ref{eq:dotTildePhi_v2}), with the initial condition $\tilde{\phi}_{\text{av}}(0) = 0$, is 
\begin{align}
\tilde{\phi}_{\text{av}}(\bar{t}) = \dfrac{\sqrt{\dfrac{p}{q}}}{1 - \dfrac{a^{2}K}{2 \omega}\max\left\{\dfrac{1}{|H^{\ast}|},1,|H^{\ast}|\right\} \dfrac{q}{p} \bar{t}} - \sqrt{\frac{p}{q}}.  \label{eq:phi-bart_v2}
\end{align}Since $\phi_{\text{av}}(t)$ in (\ref{eq:phi_1_dynamic_pf4}) is an average version of $\phi(t)=\sqrt{\frac{p}{q}}\frac{|e(\bar{t})|}{|\hat{G}(\bar{t})|}$, by invoking \cite[Theorem 2]{P:1979}, one can write  
\begin{align}
|\phi(t)-\tilde{\phi}_{\text{av}}(t)|\leq\mathcal{O}\left(\frac{1}{\omega}\right)\,.
\end{align}
By using the Triangle inequality \cite{A:1957}, one has
\begin{align}
\phi(t)&\leq\ \phi_{\text{av}}(t)+\mathcal{O}\left(\frac{1}{\omega}\right) \leq \tilde{\phi}_{\text{av}}(t)+\mathcal{O}\left(\frac{1}{\omega}\right)\! \nonumber \\
&=\dfrac{\sqrt{\frac{p}{q}}}{1 - \dfrac{a^{2}K}{2 }\max\left\{\dfrac{1}{|H^{\ast}|},1,|H^{\ast}|\right\} \frac{q}{p} t} - \sqrt{\dfrac{p}{q}} +\mathcal{O}\left(\frac{1}{\omega}\right). \label{eq:phi_t}
\end{align}
Now, defining 
\begin{align}
\hat{\phi}(t):=\dfrac{\sqrt{\dfrac{p}{q}}}{1 - \dfrac{a^{2}K}{2 }\max\left\{\dfrac{1}{|H^{\ast}|},1,|H^{\ast}|\right\} \dfrac{q}{p} t} - \sqrt{\dfrac{p}{q}} +\mathcal{O}\left(\frac{1}{\omega}\right)\,, \label{eq:hatPhi_t}
\end{align}
a lower bound on the inter-event time of original system using the static ET-ESC is given by the time it takes for the function (\ref{eq:hatPhi_t}) to go from 0 to 1.  This is at least equal to 
\begin{align} \label{cataflan_novalgina}
\tau^{\ast}&=\dfrac{2 }{a^{2}K\max\left\{\dfrac{1}{|H^{\ast}|},1,|H^{\ast}|\right\}}\frac{\beta^2}{\sigma^2}\frac{1-\mathcal{O}(1/\omega)}{1+\beta/\sigma-\mathcal{O}(1/\omega)}\,,
\end{align} 
and the Zeno behavior is avoided in the original system as well. 
\end{small}

\end{document}